\documentclass[12pt, a4paper]{amsart}

\usepackage[english]{babel}
\usepackage{amsmath}
\usepackage{amssymb}
\usepackage{amscd} 
\usepackage{microtype} 
\usepackage{verbatim} 
\usepackage{color}  
\usepackage{tikz-cd}\usetikzlibrary{babel} 
\usepackage{enumitem}
\usepackage{mathtools} 
\usepackage{cite}
\usepackage{hyperref} 
\usepackage[initials,msc-links,backrefs]{amsrefs}
\usepackage[all]{xy}

\usepackage[OT2, T1]{fontenc}

\title[On absolutely profinitely solitary lattices]
 {On absolutely profinitely solitary lattices in higher rank Lie groups}

 \author[H. Kammeyer]{Holger Kammeyer}
 
 \address{Mathematical Institute, University of D{\"u}sseldorf, Germany}
  \email{holger.kammeyer@hhu.de}
 
\subjclass[2010]{22E40, 20E18}
\keywords{profinite rigidity, lattices}

\newcounter{conj}

\theoremstyle{plain}
\newtheorem{theorem}{Theorem}

\newtheorem{conjecture}[conj]{Conjecture}
\theoremstyle{definition}

\newtheorem*{definition*}{Definition}

\newtheorem*{observation*}{Observation}

\providecommand{\ignore}[1]{}

\providecommand{\R}{\mathbb{R}}
\providecommand{\Q}{\mathbb{Q}}
\providecommand{\Z}{\mathbb{Z}}

\providecommand{\C}{\mathbb{C}}

\DeclareMathOperator{\im}{im}

\newcommand*{\arXiv}[1]{ \href{http://www.arxiv.org/abs/#1}{arXiv:\textbf{#1}}}

\begin{document}

\begin{abstract}
We establish conditions under which lattices in certain simple Lie groups are profinitely solitary in the absolute sense, so that the commensurability class of the profinite completion determines the commensurability class of the group among finitely generated residually finite groups.  While cocompact lattices are typically not absolutely solitary, we show that noncocompact lattices in \(\operatorname{Sp}(n,\R)\), \(G_{2(2)}\), \(E_8(\C)\), \(F_4(\C)\), and \(G_2(\C)\) are absolutely solitary if a well-known conjecture on Grothendieck rigidity is true.
\end{abstract}

\maketitle

\section{Introduction}

A central problem in profinite rigidity is the question whether fundamental groups of finite volume hyperbolic 3-manifolds with isomorphic profinite completions must be isomorphic themselves~\cite{Reid:profinite-rigidity}*{Section~4}.  This can be rephrased as the case \(G = \operatorname{PSL}_2(\C)\) of the general problem whether lattices \(\Gamma, \Lambda \le G\) in a connected simple Lie group \(G\) can only satisfy \(\widehat{\Gamma} \cong \widehat{\Lambda}\) if \(\Gamma \cong \Lambda\).  Recent breakthroughs are Yi Liu's result~\cite{Liu:almost} that for \(G = \operatorname{PSL}_2 (\C)\), only finitely many lattices \(\Gamma \le G\) can have the same profinite completion and M.\,Stover's negative answer~\cite{Stover:not-rigid} for \(G = \operatorname{SU}(n,1)\).  For most higher rank Lie groups~\(G\), the question can be answered in the negative by means of the \emph{congruence subgroup property} (CSP) as we recalled in~\cite{Kammeyer-Kionke:rigidity}*{Theorem 1.3}.  These non-isomorphic but profinitely isomorphic lattices are however \emph{commensurable} and hence not genuinely distinct.  So for higher rank \(G\), it appears more sincere to ask for the following, using the symbol ``\(\approx\)'' for commensurable groups.

\begin{definition*}
A lattice \(\Gamma \le G\) is called \emph{profinitely solitary} in \(G\), if every other lattice \(\Lambda \le G\) with \(\widehat{\Gamma} \approx \widehat{\Lambda}\) satisfies \(\Gamma \approx \Lambda\).
\end{definition*}

\noindent With Kionke~\citelist{\cite{Kammeyer-Kionke:rigidity} \cite{Kammeyer-Kionke:adelic}}, we found three classes of profinitely solitary lattices:
\begin{enumerate}[label=(\roman*)]
\item All lattices in simple complex Lie groups of type \(E_8\), \(F_4\), \(G_2\).
\item Noncocompact lattices in the simple real Lie groups \(E_{8(8)}\), \(E_{8(-24)}\), \(F_{4(4)}\), \(G_{2(2)}\) as well as those cocompact lattices in these groups which are defined over a Galois extension of \(\Q\).
\item Noncocompact lattices in higher rank simple real Lie groups with complexification of type \(B_n\), \(C_n\), \(E_7\).
\end{enumerate}

In this paper, we treat the naturally arising question which of these lattices \(\Gamma\) are \emph{absolutely solitary}, meaning that any finitely generated, residually finite group \(\Delta\) with \(\widehat{\Gamma} \approx \widehat{\Delta}\) satisfies \(\Gamma \approx \Delta\).  Our answer will depend on two well-known open conjectures, the first of which has appeared in~\cite{Platonov-Tavgen:grothendieck}*{Problem~2}.

\begin{conjecture} \label{conj:no-grothendieck-pairs}
  Let \(\mathbf{G}\) be an algebraic group with CSP.  Then every arithmetic subgroup \(\Lambda \le \mathbf{G}\) is \emph{Grothendieck rigid}: if an inclusion \(g \colon \Delta \hookrightarrow \Lambda\) induces an isomorphism  \(\widehat{g} \colon \widehat{\Delta} \xrightarrow{\cong} \widehat{\Lambda}\), then \(\Delta = \Lambda\).
\end{conjecture}

While our positive results are conditioned on Conjecture~\ref{conj:no-grothendieck-pairs}, the negative results are unconditional except in one case where we need to assume the following conjecture~\cite{Lubotzky:non-arithmetic}*{Section~4}.

\begin{conjecture} \label{conj:f4-csp}
  Lattices in the simple Lie group \(F_{4(-20)}\) have CSP.
\end{conjecture}

In the following theorems, all Lie groups are assumed to be connected.

\begin{theorem} \label{thm:complex-exceptional}
  Let \(\Gamma\) be a lattice in a simple complex Lie group of type \(E_8\), \(F_4\), or \(G_2\). If \(\Gamma\) is
  \begin{enumerate}[label=(\roman*),leftmargin=7mm]
  \item cocompact, then \(\Gamma\) is not absolutely solitary.
  \item noncocompact and Conjecture~\ref{conj:no-grothendieck-pairs} holds, then \(\Gamma\) is absolutely solitary.
  \end{enumerate}
\end{theorem}

\begin{theorem} \label{thm:real-exceptional}
  Let \(\Gamma\) be a lattice in a simple real Lie group of type
  \begin{enumerate}[label=(\roman*),leftmargin=8mm]
  \item \label{item:g2} \(G_{2(2)}\).  If~\(\Gamma\) is cocompact, then \(\Gamma\) is not absolutely solitary.  If~\(\Gamma\) is noncocompact and Conjecture~\ref{conj:no-grothendieck-pairs} holds, then \(\Gamma\) is absolutely solitary.
  \item \label{item:f4} \(F_{4(4)}\).  If~\(\Gamma\) is cocompact, then \(\Gamma\) is not absolutely solitary.  If~\(\Gamma\) is noncocompact and Conjecture~\ref{conj:f4-csp} holds, then \(\Gamma\) is not absolutely solitary.
  \item \label{item:e8} \(E_{8(8)}\) or \(E_{8(-24)}\).  Then \(\Gamma\) is not absolutely solitary.
  \end{enumerate}
\end{theorem}

Recall that the \(C_n\)-type in the Cartan--Killing classification of real simple Lie groups splits up into the two subtypes \(C_n\,\mathrm{I}\) given by the split form \(\operatorname{Sp}(n; \R)\) and \(C_n\,\mathrm{II}\) given by the forms \(\operatorname{Sp}(p,q)\) with \(p+q=n\).

\begin{theorem} \label{thm:other-types}
  Let \(\Gamma\) be a lattice in a higher rank simple real Lie group \(G\) with complexification of type \(B_n\), \(C_n\), or \(E_7\).
  \begin{enumerate}[label=(\roman*),leftmargin=8mm]
  \item \label{item:type-bn} In type \(B_n\), assume \(n \ge 6\) if \(\Gamma\) is cocompact and \(n \ge 14\) otherwise.
  \item \label{item:type-cn} In type \(C_n\,\mathrm{II}\), assume \(n \ge 2\) if \(\Gamma\) is cocompact and \(n \ge 6\) otherwise.
  \end{enumerate}
  Then \(\Gamma\) is not absolutely solitary.  Moreover, every cocompact lattice in \(E_{7(-5)}\) is not absolutely solitary and there exists a noncocompact lattice in a real form of type \(E_7\) which is not absolutely solitary.
\end{theorem}

The last statement is derived from the pigeonhole principle so that we cannot say in which of the real forms \(E_{7(7)}, E_{7(-5)}, E_{7(-25)}\) the lattice occurs.  In contrast, for type \(C_n\,\mathrm{I}\), we have the following result.

\begin{theorem} \label{thm:symplectic}
For \(n \ge 2\), let \(\Gamma \le \operatorname{Sp}(n;\R)\) be a noncocompact lattice and assume Conjecture~\ref{conj:no-grothendieck-pairs} holds true.  Then \(\Gamma\) is absolutely solitary.
\end{theorem}

Also recall that the Lie groups \(\operatorname{Sp}(2;\R)\) and \(\operatorname{SO}^0(3,2)\) are isogenous so that, assuming Conjecture~\ref{conj:no-grothendieck-pairs}, the case \(n=2\) of the theorem provides absolutely solitary lattices in type \(B_2\).  As another remark, the cocompact lattices in Theorem~\ref{thm:real-exceptional} which as arithmetic groups are defined over non-Galois extensions of \(\Q\), as well as all the cocompact lattices in Theorem~\ref{thm:other-types}, do not appear on the list of profinitely solitary lattices from the beginning.  While we had shown before in~\cite{Kammeyer-Kionke:rigidity}*{Theorem~1.1} that there are cocompact lattices in the corresponding Lie groups which are not profinitely solitary, it is now shown that none of the cocompact lattices are absolutely solitary.

\medskip
Let us briefly outline the proof methods.  Since \(\Gamma\) is a higher rank lattice, we may assume it is given by an arithmetic subgroup \(\Gamma \le \mathbf{H_1}\) of a simply connected almost simple linear algebraic group defined over some number field \(k\).  The key point is then a recent theorem due to R.\,Spitler~\cite{Spitler:profinite} which has emerged from his seminal joint work with M.\,Bridson--D.\,McReynolds--A.\,Reid~\cite{BMRS}.  It asserts that if \(\widehat{\Delta} \cong \widehat{\Gamma}\), then the ``group from the gutter'' \(\Delta\) embeds as a subgroup \(g \colon \Delta \hookrightarrow \Lambda\) of an arithmetic group \(\Lambda \le \mathbf{H_2}\) of an \(l\)-group \(\mathbf{H_2}\) such that we have an isomorphism of the finite adele rings \(\mathbb{A}^f_k \cong \mathbb{A}^f_l\) over which the groups are isomorphic,  \(\mathbf{H_1}(\mathbb{A}^f_k) \cong \mathbf{H_2}(\mathbb{A}^f_l)\).  If we can ensure that also \(\Lambda\) has the congruence subgroup property, only two problems remain to be addressed: Under which conditions can we conclude that \(k \cong l\) and \(\mathbf{H_1} \cong \mathbf{H_2}\) so that \(\Lambda\) and \(\Gamma\) are commensurable?  And is it true that \(\widehat{g} \colon \widehat{\Delta} \xrightarrow{\cong} \widehat{\Lambda}\) implies \(g \colon \Delta \xrightarrow{\cong} \Lambda\)?

We solve the first problem by Galois cohomological methods.  In particular, the Hasse principle and the Galois cohomology computations of simply-connected real groups in~\cite{Adams-Taibi:real} will enter repeatedly.  The second problem, Conjecture~\ref{conj:no-grothendieck-pairs}, appears however out of reach for now as we know so little about infinite index subgroups of higher rank arithmetic groups.  Note the contrast to the rank one case, where 3-manifold methods and the occurrence of nonzero first \(\ell^2\)-Betti numbers allow in some cases the construction of absolutely profinitely rigid groups~\cite{Bridson-Reid:Hopf}.

\medskip
Each of the four proofs will be given in its own section.  I wish to acknowledge financial support by the German Research Foundation, DFG 441848266 (SPP 2026/2), and DFG 284078965 (RTG 2240).  Moreover, I am grateful to the organizers of the MATRIX-MFO Tandem Workshop: ``Invariants and Structures in Low-Dimensional Topology'' in Oberwolfach in September 2021 and the members of the working group ``Profinite rigidity of 3-manifold groups'' within that workshop where the idea to this work was born:  Carol Badre, Stefan Friedl, Grace Garden, Boris Okun, Jessica Purcell, Arunima Ray, Marcy Robertson, Benjamin Ruppik, and Stephan Tillmann.

\section{Proof of Theorem~\ref{thm:complex-exceptional}}

  Let \(G\) be a connected simple complex Lie group of Cartan type \(E_8\), \(F_4\), or \(G_2\) and let \(\Gamma \le G\) be a lattice.  Note that for such a Cartan type, root lattice and weight lattice coincide and the Dynkin diagram has no symmetries.  This has the effect that the adjoint form \(\mathbf{G} = \operatorname{Aut}(\mathfrak{g})\) of the Lie algebra \(\mathfrak{g}\) of \(G\) is a connected and simply-connected \(\C\)-group with trivial center.  Hence we can identify \(G\) with \(\mathbf{G}(\C)\).  By Margulis arithmeticity~\cite{Margulis:discrete-subgroups}*{Theorem~IX.1.11, p.\,298 and (\(\ast\ast\)) as well as Remark~1.6\,(i), pp.\,293--294}, there exists a dense subfield \(k \subset \C\) with \([k : \Q] < \infty\) and there exists a connected and simply connected absolutely almost simple \(k\)-group \(\mathbf{H_1}\) such that
  \begin{enumerate}[label=(\roman*)]
  \item the complex Lie groups \(\mathbf{H_1}(\C)\) and \(\mathbf{G}(\C)\) are isomorphic,
  \item \label{item:real-places} all other infinite places of \(k\) are real and \(\mathbf{H_1}\) is anisotropic at these places,
  \item the up to commensurability defined arithmetic group \(\mathbf{H_1}(\mathcal{O}_k)\) of \(k\)-integral points of \(\mathbf{H_1}\) is commensurable with \(\Gamma\).
  \end{enumerate}

  Passing to a finite index subgroup if need be, we can thus assume~\(\Gamma\) is an arithmetic subgroup of \(\mathbf{H}_1\).  By the work of Rapinchuk~\cite{Rapinchuk:congruence}, and Raghunathan~\citelist{\cite{Raghunathan:congruence-1} \cite{Raghunathan:congruence-2}} if \(\mathbf{H_1}\) is \(k\)-isotropic, the group \(\mathbf{H_1}\) has the congruence subgroup property and it satisfies strong approximation~\cite{Platonov-Rapinchuk:algebraic-groups}*{Theorem~7.12, p.\,427}, so \(\widehat{\Gamma}\) embeds as an open subgroup of the finite adele points \(\mathbf{H_1}(\mathbb{A}_k^f)\).  More precisely, here and in the remainder the \emph{congruence subgroup property} shall mean the congruence kernel \(C(k, \mathbf{H_1})\) of \(\mathbf{H_1}\) is finite, so we can assume the canonical map \(\widehat{\Gamma} \rightarrow \mathbf{H_1}(\mathbb{A}^f_k)\) is an embedding, possibly replacing \(\Gamma\) by a finite index subgroup.

  Let \(\Delta\) be any finitely generated residually finite group whose profinite completion \(\widehat{\Delta}\) is commensurable with \(\widehat{\Gamma}\).  By~\cite{Nikolov-Segal:finitely-generated}, this is equivalent to the existence of an isomorphism from an open subgroup of \(\widehat{\Delta}\) to an open subgroup of \(\widehat{\Gamma}\).  In particular, there exist finite index subgroups \(\Delta_0 \le \Delta\) and \(\Gamma_0 \le \Gamma\) such that \(\widehat{\Delta_0} \cong \widehat{\Gamma_0}\).  By Spitler's theorem \cite{Spitler:profinite}*{Theorem~7.1}, there exists a number field \(l\), an \(l\)-form \(\mathbf{H_2}\) of \(\mathbf{G}\), an arithmetic subgroup \(\Lambda\) of \(\mathbf{H_2}\), and an injective homomorphism \(g \colon \Delta_0 \rightarrow \Lambda\) such that we have an isomorphism \(\mathbb{A}_k^f \cong \mathbb{A}_l^f\) of the finite adele rings of \(k\) and \(l\) over which \(\mathbf{H_1}(\mathbb{A}_k^f) \cong \mathbf{H_2}(\mathbb{A}_l^f)\).  Since in particular \(\mathbf{H_2}\) is isomorphic to \(\mathbf{H_1}\) over the algebraic closure \(\overline{k}\), the group \(\mathbf{H_2}\) also has type \(E_8\), \(F_4\), or \(G_2\).  Since \(k\) and \(l\) have isomorphic finite adele rings, they are in particular \emph{arithmetically equivalent}, hence have the same signature~\cite{Klingen:similarities}*{Theorem~III.1.4}.  So~\(l\) likewise has a complex place.  It follows that \(\mathbf{H_2}\) is of higher rank and therefore also has CSP.  Thus Spitler's theorem says in addition that \(g\) induces an isomorphism \(\widehat{g} \colon \widehat{\Delta_0} \xrightarrow{\cong} \widehat{\Lambda}\) on profinite completions and that \(\widehat{\Lambda} \cong \widehat{\Gamma_0}\), so \(\widehat{\Lambda}\) is commensurable with \(\widehat{\Gamma}\).

  Now assume \(\Gamma\) is noncocompact.  In that case, by~\ref{item:real-places} above, the number field \(k\) must be imaginary quadratic.  For if it had a real place, \(\mathbf{H_1}\) would be anisotropic at that place, hence anisotropic over \(k\) so \(\Gamma\) would be cocompact~\cite{Platonov-Rapinchuk:algebraic-groups}*{Theorem~4.17.(3), p.\,227}.  As number fields of degree \(\le 6\) are \emph{arithmetically solitary}~\cite{Perlis:equation}*{Section~4}, the isomorphism \(\mathbb{A}_k^f \cong \mathbb{A}_l^f\) implies \(k = l\).  So \(\Gamma\) and \(\Lambda\) are profinitely commensurable lattices in the same Lie group \(\mathbf{H_1}(\C) \cong \mathbf{H_2}(\C) \cong \mathbf{G}(\C)\), whence they are commensurable by \cite{Kammeyer-Kionke:rigidity}*{Theorem~1.2}.  Assuming Conjecture~\ref{conj:no-grothendieck-pairs}, the homomorphism \(g \colon \Delta_0 \rightarrow \Lambda\) is moreover an isomorphism.  So \(\Delta\) is commensurable with~\(\Gamma\).

  Assuming on the other hand \(\Gamma\) is cocompact, we claim that \(k\) cannot be imaginary quadratic.  Indeed, if \(k\) was imaginary quadratic, then the first noncommutative pointed Galois cohomology set \(H^1(k, \mathbf{H_1})\) would be trivial by the Hasse principle for simply-connected groups~\cite{Platonov-Rapinchuk:algebraic-groups}*{Theorem~6.6, p.\,286}.  But again the exceptional Cartan type \(E_8\), \(F_4\), or \(G_2\) of \(\mathbf{H}_1\) shows that \(\operatorname{Aut}_{\overline{k}}(\mathbf{H_1}) \cong \mathbf{H_1}\).  Since \(H^1(k, \operatorname{Aut}_{\overline{k}}(\mathbf{H_1}))\) classifies the \(k\)-twists of \(\mathbf{H_1}\), the triviality of this set would show that \(\mathbf{H}_1\) was \(k\)-split, hence \(\Gamma\) would not be cocompact.  We thus have \([k : \Q] \ge 3\).  By~\ref{item:real-places} above, \(k\) must have precisely one complex place and at least one real place and \(\mathbf{H_1}\) is anisotropic at all real places.  Let \(\mathbf{H_2}\) be the \(k\)-split form of \(\mathbf{G}\) which up to commensurability defines the arithmetic subgroup \(\Lambda = \mathbf{H_2}(\mathcal{O}_k)\).  Since both \(\mathbf{H_1}\) and \(\mathbf{H_2}\) split at all finite places of \(k\) by the Galois cohomology computation of Kneser~\cite{Kneser:galois-2}*{table on p.\,254}, we conclude \(\mathbf{H_1}(\mathbb{A}_k^f) \cong \mathbf{H_2}(\mathbb{A}_k^f)\) from \cite{Kammeyer-Kionke:adelic}*{Lemma~2.5 and Lemma~2.2}.  This, the congruence subgroup property, and strong approximation imply that \(\widehat{\Lambda}\) is commensurable with \(\widehat{\Gamma}\).  If \(\Gamma\) and \(\Lambda\) were themselves commensurable, however, strong rigidity~\cite{Margulis:discrete-subgroups}*{Theorem~7.1, p.\,251} would imply that \(\mathbf{H_1}\) and \(\mathbf{H_2}\) split at the same number of real places contradicting that \(\mathbf{H_1}\) is anisotropic at all real places while \(\mathbf{H_2}\) splits at all real places.

  \section{Proof of Theorem~\ref{thm:real-exceptional}}
    
  Let \(\Gamma \le G\) be a lattice in a connected simple real Lie group \(G\) of either of the four listed types.  Replacing \(\Gamma\) by a finite index subgroup if need be, we may assume that \(G\) is the adjoint form (has trivial center).  As the Dynkin diagrams of types \(E_8\), \(F_4\), and \(G_2\) have no symmetries, \(G\) can be realized as \(G = \mathbf{G}(\R)\) where \(\mathbf{G} = \operatorname{Aut}_\R(\mathfrak{g})\) and \(\mathfrak{g}\) denotes the Lie algebra of \(G\).  Margulis arithmeticity~\cite{Margulis:discrete-subgroups}*{(\(\ast\ast\)) and Remark~1.6\,(i), pp.\,293--294} now says that there exists a totally real number field \(k\) with a distinguished embedding \(k \subset \R\) and a simply-connected \(k\)-group \(\mathbf{H}_1\) such that \(\mathbf{H}_1(\R)\) is isogenous to \(\mathbf{G}(\R)\), such that \(\mathbf{H_1}\) is anisotropic at all the other infinite places, and such that \(\mathbf{H_1}(\mathcal{O}_k)\) is commensurable with \(\Gamma\).  By~\cite{Platonov-Rapinchuk:algebraic-groups}*{Theorems~9.23 and~9.24, pp.\,568--569}, which were proven by Rapinchuk~\citelist{\cite{Rapinchuk:affine} \cite{Rapinchuk:congruence}} in the exceptional case, the group~\(\mathbf{H_1}\) has the congruence subgroup property, so \(\widehat{\Gamma}\) is an open subgroup of~\(\mathbf{H_1}(\mathbb{A}^f_k)\).

    Now we assume \(\Gamma\) is cocompact.  We claim that then \([k : \Q] \ge 2\).  Indeed, since \(\mathbf{H}_1\) is a \(k\)-form of type \(E_8\), \(F_4\), or \(G_2\), we have \(\operatorname{Aut}_{\overline{k}}(\mathbf{H_1}) \cong_{\,\overline{k}} \mathbf{H_1}\).  So if we had \(k = \Q\), then the Hasse principle~\cite{Platonov-Rapinchuk:algebraic-groups}*{Theorem~6.6, p.\,286} for Galois cohomology of simply-connected groups would show \(H^1(\Q, \mathbf{H_1}) \cong H^1(\R, \mathbf{H_1})\) so that the \(\Q\)-forms of \(\mathbf{H_1}\) are in one to one correspondence with the \(\R\)-forms.  It follows from a criterion of Prasad--Rapinchuk~\cite{Prasad-Rapinchuk:isotropic}*{Theorem~1\,(iii)}, that this correspondence takes \(\R\)-isotropic forms to \(\Q\)-isotropic forms.  Therefore \(\mathbf{H_1}\) would be \(\Q\)-isotropic and \(\Gamma\) would not be cocompact.  We thus have another real place \(v\) of \(k\) apart from the distinguished one.  Then as in the proof of Theorem~\ref{thm:complex-exceptional}, the \(k\)-split form \(\mathbf{H_2}\) of \(\mathbf{G}\) gives an arithmetic subgroup \(\Lambda = \mathbf{H_2}(\mathcal{O}_k)\) such that \(\widehat{\Lambda}\) is commensurable with \(\widehat{\Gamma}\), but \(\Gamma\) is not commensurable with \(\Lambda\) by strong rigidity because \(\mathbf{H_1}\) is anisotropic at \(v\) while \(\mathbf{H_2}\) splits at \(v\).

    Assuming \(\Gamma\) is noncocompact, \(\mathbf{H_1}\) must be \(k\)-isotropic, hence this time the condition that \(\mathbf{H_1}\) is anisotropic at each non-distinguished infinite place implies~\(k = \Q\).  Let \(\Delta\) be a finitely generated residually finite group such that \(\widehat{\Delta}\) is commensurable with \(\widehat{\Gamma}\).  Reasoning as in the proof of Theorem~\ref{thm:complex-exceptional}, we obtain a \(\Q\)-form \(\mathbf{H}_2\) of \(\mathbf{G}\), an arithmetic subgroup \(\Lambda\) of \(\mathbf{H_2}\), an embedding \(g \colon \Delta_0 \rightarrow \Lambda\) of a finite index subgroup \(\Delta_0 \le \Delta\), and an isomorphism \(\mathbf{H_1}(\mathbb{A}^f_\Q) \cong \mathbf{H_2}(\mathbb{A}^f_\Q)\).  In particular \(\mathbf{H_2}\) has the same complex Cartan type as \(\mathbf{H_1}\).

    Let us now first assume that \(\mathbf{H_1}\) and hence \(\mathbf{H_2}\) have type \(G_2\).  We have \(H^1(\Q, \operatorname{Aut}_{\overline{\Q}}(\mathbf{H_1})) \cong H^1(\Q, \mathbf{H_1}) \cong H^1(\R, \mathbf{H_1})\) and the \(\R\)-split form is the only isotropic \(\R\)-form in type \(G_2\).  Since both \(\mathbf{H_1}\) and \(\mathbf{H_2}\) are \(\R\)-isotropic, it follows that \(\mathbf{H_1}\) is \(\Q\)-isomorphic to \(\mathbf{H_2}\), whence \(\Gamma\) and \(\Lambda\) are commensurable and \(\mathbf{H_2}\) has CSP.   So \(g\) induces an isomorphism \(\widehat{\Delta_0} \cong \widehat{\Lambda}\) and if Conjecture~\ref{conj:no-grothendieck-pairs} holds true, we conclude that \(\Delta\) is commensurable with \(\Gamma\).

    Now we assume that \(\mathbf{H_1}\) is of type \(E_8\) or \(F_4\).  Then there exists an \(\R\)-isotropic form different from \(\mathbf{H_1}(\R)\) and hence a simply-connected \(\Q\)-form \(\mathbf{H_2}\) with \(\mathbf{H_2}(\R) \not\cong \mathbf{H_1}(\R)\).  To construct \(\mathbf{H_2}\), one can for example take the simply-connected covering group of the automorphism group of the \(\Q\)-Lie algebra of the other \(\R\)-form as constructed in~\cite{Borel:clifford}*{Proposition~3.7}.  We have \(\mathbf{H_1}(\mathbb{A}^f_\Q) \cong \mathbf{H_2}(\mathbb{A}^f_\Q)\) again by \cite{Kammeyer-Kionke:adelic}*{Lemma~2.5 and Lemma~2.2} because \(p\)-adic forms in type \(E_8\) and \(F_4\) are unique by Kneser's computation~\cite{Kneser:galois-2}*{table on p.\,254}.  In type \(F_4\), the real group \(\mathbf{H_2}(\R)\) is the rank one group \(F_{4(-20)}\) so that we need to assume Conjecture~\ref{conj:f4-csp} to ensure \(\mathbf{H_2}\) has the congruence subgroup property.  In type \(E_8\), this is granted because \(\mathbf{H_2}(\R)\) has higher rank.  Let \(\Lambda\) be the up to commensurability defined group \(\mathbf{H_2}(\Z)\).  The congruence subgroup property effects that \(\widehat{\Gamma}\) and \(\widehat{\Lambda}\) are open subgroups of \(\mathbf{H_1}(\mathbb{A}^f_\Q)\) and \(\mathbf{H_2}(\mathbb{A}^f_\Q)\), respectively, hence these profinite completions are commensurable.  But the groups \(\Gamma\) and \(\Lambda\) are not commensurable themselves by strong rigidity, being lattices in different simple Lie groups.

    \section{Proof of Theorem~\ref{thm:other-types}}
    
  Once again, we may assume \(\Gamma \le \mathbf{G}(\R)\) where \(\mathbf{G} = \operatorname{Aut}_\R(\mathfrak{g})\) and \(\mathfrak{g}\) is the Lie algebra of \(G\).  Margulis arithmeticity gives a totally real number field \(k\) with a distinguished embedding \(k \subset \R\) and a simply-connected \(k\)-group \(\mathbf{H_1}\) such that \(\mathbf{H_1}(\R)\) covers \(\mathbf{G}(\R)\), such that \(\mathbf{H_1}\) is anisotropic at all the other real places of \(k\), and such that \(\mathbf{H_1}(\mathcal{O}_k)\) is commensurable with \(\Gamma\).  By~\cite{Platonov-Rapinchuk:algebraic-groups}*{Theorems~9.23 and~9.24, pp.\,568--569}, \(\mathbf{H_1}\) has the congruence subgroup property and hence \(\widehat{\Gamma}\) is an open subgroup of \(\mathbf{H_1}(\mathbb{A}^f_k)\).

  The \(k\)-twists of \(\mathbf{H_1}\) are classified by the first pointed Galois cohomology set \(H^1(k, \operatorname{Aut}(\mathbf{H_1}))\).  Here \(\operatorname{Aut}(\mathbf{H}_1) \cong \overline{\mathbf{H_1}} = \mathbf{H_1} / \mathbf{Z}\) can be identified with the adjoint form \(\overline{\mathbf{H_1}}\) where \(\mathbf{Z} \subset \mathbf{H_1}\) denotes the center because the Dynkin diagrams of type \(B_n\), \(C_n\), and \(E_7\) have no symmetries.  By~\cite{Serre:galois}*{Proposition~43, p.\,55}, for every field extension \(K/k\), we have an exact sequence of pointed sets
    \[ H^1(K, \mathbf{Z}) \longrightarrow H^1(K, \mathbf{H_1}) \xrightarrow{\ \pi_K \ } H^1(K, \overline{\mathbf{H_1}}) \xrightarrow{\ \delta_K \ } H^2(K, \mathbf{Z}). \]
    For the distinguished extension \(\R / k\), we claim that \(\ker \delta_\R\) consists of at least four elements in the cases listed in~\ref{item:type-bn} and~\ref{item:type-cn} of the Theorem.  By~\cite{Serre:galois}*{Proposition~42, p.\,54}, the number \(|\ker \delta_K| = |\im \pi_K|\) equals the number of \(H^1(K,\mathbf{Z})\)-orbits in \(H^1(K, \mathbf{H_1})\).  For the Cartan type \(B_n\), \(C_n\), or \(E_7\) of \(\mathbf{H_1}\), we have \(\mathbf{Z} = \{\pm 1\}\), see for instance~\cite{Platonov-Rapinchuk:algebraic-groups}*{table on p.\,332}, hence also \(H^1(\R, \mathbf{Z}) = \{ \pm 1 \}\) by~\cite{Platonov-Rapinchuk:algebraic-groups}*{Lemma~2.6, p.\,73}.  The real Galois cohomology of simply-connected \(\R\)-groups was first computed by Borovoi~\cite{Borovoi:real}; tables can be found in the recent work of Adams--Ta{\"i}bi~\cite{Adams-Taibi:real}*{Section~10, Tables~1 and~2}.  We have \(|H^1(\R, \mathbf{H_1})| \ge \lfloor \frac{n}{2} \rfloor\) if \(\mathbf{H_1}(\R) \cong \operatorname{Spin}(p,q)\) with \(p+q = 2n+1\) and \(|H^1(\R, \mathbf{H_1})| = n+1\) if \(\mathbf{H_1}(\R) \cong \operatorname{Sp}(p,q)\) with \(p+q = n\).  Since \(H^1(\R, \mathbf{Z})\)-orbits have at most two elements, we see that for \(n \ge 14\) in the first case and for \(n \ge 6\) in the second case, we have \(|\ker \delta_\R| \ge 4\). 

    Now since \(\ker \delta_\R\) has at least four elements, there exists a nontrivial class \([a_\infty] \in \ker \delta_\R\) which neither corresponds to the unique compact real form nor to the unique real form of real rank one.  Fix any finite place \(v_0\) of \(k\).  We have a commutative square
  \[
    \begin{tikzcd}
  H^1(k, \overline{\mathbf{H_1}}) \ar[r, "\alpha"] \ar[d, "\delta_k"] & \bigoplus\limits_{v \neq v_0} H^1(k_v, \overline{\mathbf{H_1}}) \ar[d, "\bigoplus\limits_{v \neq v_0} \delta_{k_v}"] \\
  H^2(k, \mathbf{Z}) \ar[r, "\beta"] & \bigoplus\limits_{v \neq v_0} H^2(k_v, \mathbf{Z})
    \end{tikzcd}
  \]
  and by~\cite{Prasad-Rapinchuk:isotropic}*{Proposition~1 and Theorem~3}, \(\alpha\) is surjective and both \(\alpha\) and \(\beta\) are injective.  Hence there exists a unique global class \([a] \in H^1(k, \overline{\mathbf{H_1}})\) such that \([a]\) localizes to \([a_\infty] \in H^1(\R, \overline{\mathbf{H_1}})\) at the distinguished real place and to the unit class \([1] \in H^1(k_v, \overline{\mathbf{H_1}})\) at all other places \(v \neq v_0\) of \(k\).  We then have
  \[ 1 = \left(\textstyle \bigoplus_{v \neq v_0} \delta_{k_v} \right) (\alpha ([a])) = \beta (\delta_k ([a])) \]
  so \(\delta_k([a]) = 1\) which gives in particular \(\delta_{k_{v_0}} ([a]) = 1\).  As \(v_0\) is a finite place, we have M.\,Kneser's result \(H^1(k_{v_0},\mathbf{H_1}) = \{1\}\), see for instance~\cite{Platonov-Rapinchuk:algebraic-groups}*{Theorem~6.4, p.\,284}, so \(\delta_{k_{v_0}}\) is injective and \(\delta_{k_{v_0}} ([a]) = 1\) implies that \([a]\) also localizes to the unit class in \(H^1(k_{v_0}, \overline{\mathbf{H_1}})\).  This shows that the simply-connected \(k\)-form \(\mathbf{H_2}\) corresponding to \([a] \in H^1(k, \overline{\mathbf{H_1}})\) satisfies \(\mathbf{H_1}(k_v) \cong \mathbf{H_2}(k_v)\) for all finite places \(v\) of \(k\) and hence \(\mathbf{H_1}(\mathbb{A}^f_k) \cong \mathbf{H_2}(\mathbb{A}^f_k)\) by \cite{Kammeyer-Kionke:adelic}*{Lemma~2.5 and Lemma~2.2}.  Since at the distinguished real place, we have \(\operatorname{rank}_\R \mathbf{H_2}(\R) \ge 2\) and \(\mathbf{H_2}\) is of type \(B_n\) or \(C_n\), the group \(\mathbf{H_2}\) has the congruence subgroup property so for \(\Lambda = \mathbf{H_2}(\mathcal{O}_k)\), the profinite completion \(\widehat{\Lambda}\) is an open subgroup of \(\mathbf{H_2}(\mathbb{A}^f_k)\).  We conclude that \(\widehat{\Gamma}\) and \(\widehat{\Lambda}\) are commensurable but again, \(\Gamma\) and \(\Lambda\) are lattices in different real Lie groups whence not commensurable by strong rigidity.

  If we assume \(\Gamma\) is cocompact, we can improve the construction as follows.  First we show again that then of necessity \([k : \Q] \ge 2\).  Indeed, we infer from~\cite{Tits:classification}*{Table~II} that in type \(B_n\) with \(n \ge 2\), the first node in every Tits index over \(\Q_p\) is circled wheres in type \(C_n\) with \(n \ge 3\) and type \(E_7\), the second node in every Tits index over \(\Q_p\) is circled.  Since the same is true for every isotropic Tits index of these types over \(\R\), a criterion of Prasad--Rapinchuk~\cite{Prasad-Rapinchuk:isotropic}*{Theorem~1\,(ii) and\,(iii)} and the injectivity of the map \(\alpha\) above show that every \(\R\)-isotropic \(\Q\)-form is actually \(\Q\)-isotropic, hence defines a noncocompact arithmetic group.

  Thus we have another real place \(v_1\) of \(k\) at which \(\mathbf{H_1}\) is anisotropic.  We now claim that \(\ker \delta_{k_{v_1}}\) is not trivial.  Indeed, we have \(\mathbf{H_1}(k_{v_1}) \cong \operatorname{Spin}(2n+1)\) in type \(B_n\) and \(\mathbf{H_1}(k_{v_1}) \cong \operatorname{Sp}(n)\) in type \(C_n\).  Thus by the above consideration, \(\ker \delta_{k_{v_1}}\) is definitely nontrivial if \(n \ge 6\) or \(n \ge 2\), respectively.  Moreover we have \(|H^1(\R; \mathbf{H_1})| = 4\) if \(\mathbf{H_1}(\R)\) has type \(E_{7(-5)}\) by~\cite{Adams-Taibi:real}*{Section~10, Table~3}, hence \(\ker \delta_{k_{v_1}}\) is also nontrivial in that case.  Arguing as above, we obtain a \(k\)-group \(\mathbf{H_2}\) isomorphic to \(\mathbf{H_1}\) at all places except \(v_1\) where \(\mathbf{H_2}\) is a noncompact Lie group.  It follows that for \(\Lambda = \mathbf{H_2}(\mathcal{O}_k)\), the profinite completion \(\widehat{\Lambda}\) is commensurable with \(\widehat{\Gamma}\) but \(\Gamma\) and \(\Lambda\) are not commensurable, being lattices in Lie groups with different real rank.

  Finally, we show that there exists a real form of type \(E_7\) with a noncocompact lattice \(\Gamma\) which is not absolutely solitary.  Let \(\mathbf{H_1}\) be a \(\Q\)-split form of type \(E_7\) defining (up to commensurability) the arithmetic subgroup \(\Gamma = \mathbf{H_1}(\mathbb{Z})\).  Then for any fixed prime \(l\), we now have a bijection
    \[  H^1(\Q, \overline{\mathbf{H_1}}) \xrightarrow{\ \alpha \ } \bigoplus_{p \neq l} H^1(\Q_p, \overline{\mathbf{H_1}}). \]
We note from~\cite{Tits:classification}*{Table~II} that in type \(E_7\), there exist two higher rank non-split real forms but only one non-split \(p\)-adic form.  Let \([a_1], [a_2] \in H^1(\Q, \overline{\mathbf{H_1}})\) be the unique \(\Q\)-twists which map via \(\alpha\) to the trivial twist at all \(\Q_p\) for \(p \neq l\) and to the two non-split isotropic real forms at the infinite place.  Then either \([a_1]\) or \([a_2]\) also localizes to the trivial twist at \(\Q_l\) or both localize to the non-split form at \(\Q_l\).  In the former case, we obtain a \(\Q\)-form \(\mathbf{H_2}\) defining \(\Lambda = \mathbf{H_2}\) such that \(\widehat{\Lambda}\) and \(\widehat{\Gamma}\) are commensurable but \(\Gamma\) and \(\Lambda\) are not commensurable being lattices in different Lie groups.  In the latter case the \(\Q\)-forms corresponding to \([a_1]\) and \([a_2]\) define profinitely commensurable lattices which are not commensurable.

\section{Proof of Theorem~\ref{thm:symplectic}}

  Being a noncocompact lattice, by Margulis arithmeticity \(\Gamma\) is commensurable with \(\mathbf{H_1}(\Z)\) for a simply-connected \(\Q\)-group \(\mathbf{H_1}\) such that \(\mathbf{H_1}(\R) \cong \operatorname{Sp}(n, \R)\).  Since \(n \ge 2\), the group \(\mathbf{H_1}\) has the congruence subgroup property, so \(\widehat{\Gamma}\) is an open subgroup of \(\mathbf{H_1}(\mathbb{A}^f_\Q)\).  As in the proof of Theorem~\ref{thm:real-exceptional}, given a finitely generated residually finite group \(\Delta\) such that \(\widehat{\Delta}\) is commensurable with \(\widehat{\Gamma}\), we obtain a \(\Q\)-group \(\mathbf{H_2}\), an arithmetic subgroup \(\Lambda\) of \(\mathbf{H_2}\), an embedding \(g \colon \Delta_0 \rightarrow \Lambda\) of a subgroup \(\Delta_0 \le \Delta\) of finite index, and finally an isomorphism \(\mathbf{H_1}(\mathbb{A}^f_\Q) \cong \mathbf{H_2}(\mathbb{A}^f_\Q)\) showing that \(\mathbf{H_2}\) is of type \(C_n\) as well.  More precisely, the latter isomorphism shows that the class \([a] \in H^1(\Q, \overline{\mathbf{H_1}})\) defined by the cocycle \(a\) corresponding to \(\mathbf{H_2}\) lies in the kernel of the map \(\alpha\) in the diagram of pointed Galois cohomology sets
  \[
\begin{tikzcd}
{H^1(\mathbb{Q}, \mathbf{H_1})} \arrow[d] \arrow[r] & {H^1(\mathbb{Q}, \overline{\mathbf{H_1}})} \arrow[d, "\alpha"] \arrow[r, "\delta"] & {H^2(\mathbb{Q}, \mathbf{Z})} \arrow[d, "\beta", "\cong"'] \\
{\bigoplus\limits_{p \nmid \infty} H^1(\mathbb{Q}_p, \mathbf{H_1})} \arrow[r] & {\bigoplus\limits_{p \nmid \infty} H^1(\mathbb{Q}_p, \overline{\mathbf{H_1}})} \arrow[r, "\cong"', "\oplus_p\, \delta_p"] & {\bigoplus\limits_{p \nmid \infty} H^2(\mathbb{Q}_p, \mathbf{Z})}
\end{tikzcd}
\]
with exact rows.  The map \(\beta\) is a bijection as can be inferred from the Hasse--Brauer--Noether theorem~\cite{Prasad-Rapinchuk:isotropic}*{Lemma~2}.  Each of the maps \(\delta_p \colon H^1(\mathbb{Q}_p, \overline{\mathbf{H_1}}) \xrightarrow{\cong} H^2(\mathbb{Q}_p, \mathbf{Z})\) is a bijection according to~\cite{Platonov-Rapinchuk:algebraic-groups}*{Corollary to Theorem~6.20, p.\,326}.  It follows that \(\ker \alpha = \ker \delta\).  As above, the Hasse principle for simply-connected groups~\cite{Platonov-Rapinchuk:algebraic-groups}*{Theorem~6.6, p.\,286} gives \(H^1(\mathbb{Q}, \mathbf{H_1}) \cong H^1(\mathbb{R}, \mathbf{H_1})\).  But the latter set is trivial by \cite{Adams-Taibi:real}*{Section~10, Table~1} because \(\mathbf{H_1}(\R) = \operatorname{Sp}(n, \R)\).  By exactness, \(\ker \delta\) and hence \(\ker \alpha\) is trivial.  It follows that \([a]\) is the unit class, so \(\mathbf{H_2}\) is \(\Q\)-isomorphic to \(\mathbf{H_1}\), hence \(\Gamma\) is commensurable with \(\Lambda\) and \(\mathbf{H_2}\) has CSP.  So \(g\) induces an isomorphism \(\widehat{\Delta_0} \xrightarrow{\cong} \widehat{\Lambda}\) and if Conjecture~\ref{conj:no-grothendieck-pairs} is true, \(g\) is itself an isomorphism, so \(\Delta\) is commensurable with \(\Gamma\).


\begin{bibdiv}[References]

  \begin{biblist}

    \bib{Adams-Taibi:real}{article}{
   author={Adams, Jeffrey},
   author={Ta\"{\i}bi, Olivier},
   title={Galois and Cartan cohomology of real groups},
   journal={Duke Math. J.},
   volume={167},
   date={2018},
   number={6},
   pages={1057--1097},
   issn={0012-7094},
   review={\MR{3786301}},
   doi={10.1215/00127094-2017-0052},
 }
 
\bib{Borel:clifford}{article}{
   author={Borel, Armand},
   title={Compact Clifford-Klein forms of symmetric spaces},
   journal={Topology},
   volume={2},
   date={1963},
   pages={111--122},
   issn={0040-9383},
   review={\MR{146301}},
   doi={10.1016/0040-9383(63)90026-0},
 }

 \bib{Borovoi:real}{article}{
   author={Borovo\u{\i}, M. V.},
   title={Galois cohomology of real reductive groups and real forms of
   simple Lie algebras},
   language={Russian},
   journal={Funktsional. Anal. i Prilozhen.},
   volume={22},
   date={1988},
   number={2},
   pages={63--64},
   issn={0374-1990},
   translation={
      journal={Funct. Anal. Appl.},
      volume={22},
      date={1988},
      number={2},
      pages={135--136},
      issn={0016-2663},
   },
   review={\MR{947609}},
   doi={10.1007/BF01077606},
 }
 
\bib{BMRS}{article}{
   author={Bridson, M. R.},
   author={McReynolds, D. B.},
   author={Reid, A. W.},
   author={Spitler, R.},
   title={Absolute profinite rigidity and hyperbolic geometry},
   journal={Ann. of Math. (2)},
   volume={192},
   date={2020},
   number={3},
   pages={679--719},
   issn={0003-486X},
   review={\MR{4172619}},
   doi={10.4007/annals.2020.192.3.1},
}

\bib{Bridson-Reid:Hopf}{article}{
  author={Bridson, M. R.},
  author={Reid, A. W.},
  title={Profinite rigidity, Kleinian groups, and the cofinite Hopf property},
  journal={To appear in a special issue of the Michigan Math. J. honoring Gopal Prasad},
  year={2021},
  review={\arXiv{2107.14696}},
}

\bib{Kammeyer-Kionke:adelic}{article}{
   author={Kammeyer, Holger},
   author={Kionke, Steffen},
   title={Adelic superrigidity and profinitely solitary lattices},
   journal={Pacific J. Math.},
   volume={313},
   date={2021},
   number={1},
   pages={137--158},
   issn={0030-8730},
   review={\MR{4313430}},
   doi={10.2140/pjm.2021.313.137},
}

\bib{Kammeyer-Kionke:rigidity}{article}{
  author={Kammeyer, Holger},
  author={Kionke, Steffen},
  title={On the profinite rigidity of lattices in higher rank Lie groups},
  year={2020},
  review={\arXiv{2009.13442}},
}

\bib{Klingen:similarities}{book}{
   author={Klingen, Norbert},
   title={Arithmetical similarities},
   series={Oxford Mathematical Monographs},
   note={Prime decomposition and finite group theory;
   Oxford Science Publications},
   publisher={The Clarendon Press, Oxford University Press, New York},
   date={1998},
   pages={x+275},
   isbn={0-19-853598-8},
   review={\MR{1638821}},
}

\bib{Kneser:galois-2}{article}{
   author={Kneser, Martin},
   title={Galois-Kohomologie halbeinfacher algebraischer Gruppen \"{u}ber
   ${\germ p}$-adischen K\"{o}rpern. II},
   language={German},
   journal={Math. Z.},
   volume={89},
   date={1965},
   pages={250--272},
   issn={0025-5874},
   review={\MR{188219}},
   doi={10.1007/BF02116869},
 }

 \bib{Liu:almost}{article}{
  author={Liu, Yi},
  title={Finite-volume hyperbolic 3-manifolds are almost determined by their finite quotient groups},
  year={2020},
  review={\arXiv{2011.09412}},
}

\bib{Lubotzky:non-arithmetic}{article}{
   author={Lubotzky, Alexander},
   title={Some more non-arithmetic rigid groups},
   conference={
      title={Geometry, spectral theory, groups, and dynamics},
   },
   book={
      series={Contemp. Math.},
      volume={387},
      publisher={Amer. Math. Soc., Providence, RI},
   },
   date={2005},
   pages={237--244},
   review={\MR{2180210}},
   doi={10.1090/conm/387/07244},
}
 
\bib{Margulis:discrete-subgroups}{book}{
   author={Margulis, G. A.},
   title={Discrete subgroups of semisimple Lie groups},
   series={Ergebnisse der Mathematik und ihrer Grenzgebiete (3) [Results in
   Mathematics and Related Areas (3)]},
   volume={17},
   publisher={Springer-Verlag, Berlin},
   date={1991},
   pages={x+388},
   isbn={3-540-12179-X},
   review={\MR{1090825}},
   doi={10.1007/978-3-642-51445-6},
}

\bib{Nikolov-Segal:finitely-generated}{article}{
   author={Nikolov, Nikolay},
   author={Segal, Dan},
   title={On finitely generated profinite groups. I. Strong completeness and
   uniform bounds},
   journal={Ann. of Math. (2)},
   volume={165},
   date={2007},
   number={1},
   pages={171--238},
   issn={0003-486X},
   review={\MR{2276769}},
   doi={10.4007/annals.2007.165.171},
 }
 
\bib{Platonov-Tavgen:grothendieck}{article}{
   author={Platonov, V. P.},
   author={Tavgen\cprime , O. I.},
   title={Grothendieck's problem on profinite completions and
   representations of groups},
   journal={$K$-Theory},
   volume={4},
   date={1990},
   number={1},
   pages={89--101},
   issn={0920-3036},
   review={\MR{1076526}},
   doi={10.1007/BF00534194},
}

\bib{Perlis:equation}{article}{
   author={Perlis, Robert},
   title={On the equation $\zeta _{K}(s)=\zeta _{K'}(s)$},
   journal={J. Number Theory},
   volume={9},
   date={1977},
   number={3},
   pages={342--360},
   issn={0022-314X},
   review={\MR{447188}},
   doi={10.1016/0022-314X(77)90070-1},
 }

\bib{Platonov-Rapinchuk:algebraic-groups}{book}{
   author={Platonov, Vladimir},
   author={Rapinchuk, Andrei},
   title={Algebraic groups and number theory},
   series={Pure and Applied Mathematics},
   volume={139},
   note={Translated from the 1991 Russian original by Rachel Rowen},
   publisher={Academic Press, Inc., Boston, MA},
   date={1994},
   pages={xii+614},
   isbn={0-12-558180-7},
   review={\MR{1278263}},
}
 
 \bib{Prasad-Rapinchuk:isotropic}{article}{
   author={Prasad, Gopal},
   author={Rapinchuk, Andrei S.},
   title={On the existence of isotropic forms of semi-simple algebraic
   groups over number fields with prescribed local behavior},
   journal={Adv. Math.},
   volume={207},
   date={2006},
   number={2},
   pages={646--660},
   issn={0001-8708},
   review={\MR{2271021}},
   doi={10.1016/j.aim.2006.01.001},
}

\bib{Raghunathan:congruence-1}{article}{
   author={Raghunathan, M. S.},
   title={On the congruence subgroup problem},
   journal={Inst. Hautes \'{E}tudes Sci. Publ. Math.},
   number={46},
   date={1976},
   pages={107--161},
   issn={0073-8301},
   review={\MR{507030}},
 }

 \bib{Raghunathan:congruence-2}{article}{
   author={Raghunathan, M. S.},
   title={On the congruence subgroup problem. II},
   journal={Invent. Math.},
   volume={85},
   date={1986},
   number={1},
   pages={73--117},
   issn={0020-9910},
   review={\MR{842049}},
   doi={10.1007/BF01388793},
}

\bib{Rapinchuk:congruence}{article}{
   author={Rapinchuk, A. S.},
   title={On the congruence subgroup problem for algebraic groups},
   language={Russian},
   journal={Dokl. Akad. Nauk SSSR},
   volume={306},
   date={1989},
   number={6},
   pages={1304--1307},
   issn={0002-3264},
   translation={
      journal={Soviet Math. Dokl.},
      volume={39},
      date={1989},
      number={3},
      pages={618--621},
      issn={0197-6788},
   },
   review={\MR{1015345}},
 }

 \bib{Rapinchuk:affine}{article}{
   author={Rapinchuk, A. S.},
   title={The congruence subgroup problem for algebraic groups and strong
   approximation in affine varieties},
   language={Russian, with English summary},
   journal={Dokl. Akad. Nauk BSSR},
   volume={32},
   date={1988},
   number={7},
   pages={581--584, 668},
   issn={0002-354X},
   review={\MR{954097}},
}

\bib{Reid:profinite-rigidity}{article}{
   author={Reid, Alan W.},
   title={Profinite rigidity},
   conference={
      title={Proceedings of the International Congress of
      Mathematicians---Rio de Janeiro 2018. Vol. II. Invited lectures},
   },
   book={
      publisher={World Sci. Publ., Hackensack, NJ},
   },
   date={2018},
   pages={1193--1216},
   review={\MR{3966805}},
}

\bib{Serre:galois}{book}{
   author={Serre, Jean-Pierre},
   title={Galois cohomology},
   series={Springer Monographs in Mathematics},
   edition={Corrected reprint of the 1997 English edition},
   note={Translated from the French by Patrick Ion and revised by the
   author},
   publisher={Springer-Verlag, Berlin},
   date={2002},
   pages={x+210},
   isbn={3-540-42192-0},
   review={\MR{1867431}},
 }
 
\bib{Spitler:profinite}{thesis}{
author = {Ryan F. Spitler},
title = {Profinite Completions and Representations of Finitely Generated Groups},
year = {2019},
note = {PhD thesis},
organization = {Purdue University},
doi = {10.25394/PGS.9117068.v1}
}

\bib{Stover:not-rigid}{article}{
   author={Stover, Matthew},
   title={Lattices in ${\rm PU}(n,1)$ that are not profinitely rigid},
   journal={Proc. Amer. Math. Soc.},
   volume={147},
   date={2019},
   number={12},
   pages={5055--5062},
   issn={0002-9939},
   review={\MR{4021068}},
   doi={10.1090/proc/14763},
 }
 
\bib{Tits:classification}{article}{
   author={Tits, J.},
   title={Classification of algebraic semisimple groups},
   conference={
      title={Algebraic Groups and Discontinuous Subgroups},
      address={Proc. Sympos. Pure Math., Boulder, Colo.},
      date={1965},
   },
   book={
      publisher={Amer. Math. Soc., Providence, R.I., 1966},
   },
   date={1966},
   pages={33--62},
   review={\MR{0224710}},
 }
 
\end{biblist}
\end{bibdiv}
\end{document}